\documentclass[11pt,a4paper]{article}
\usepackage[english]{babel}
\usepackage[utf8]{inputenc}

\usepackage{amsmath}
\usepackage{amssymb}
\usepackage{amsfonts}
\usepackage{amsbsy}

\usepackage{graphicx}

\usepackage{newlfont}
\usepackage{float}
\usepackage[paperwidth=192mm,
            paperheight=262mm,
            vmargin={19mm,19mm},
            hmargin={13.7mm,13.7mm},
            headsep=12pt,
            footskip=12pt]{geometry}

\usepackage{hyperref}
\hypersetup{colorlinks,
            linkcolor=blue,
            citecolor=blue,
            urlcolor=blue}

\usepackage{siunitx}
\usepackage{subcaption}

\usepackage{caption}
\usepackage{multirow}

\usepackage{booktabs}
\usepackage{tabularx}

\usepackage{placeins}

\usepackage[normalem]{ulem}

\usepackage[backend=bibtex,style=numeric,sorting=nyt]{biblatex}
\renewbibmacro{in:}{}
\DeclareMathOperator{\grad}{\nabla}
\DeclareMathOperator{\dive}{\nabla\cdot}
\DeclareMathOperator{\vel}{\mathbf{u}}

\newcommand{\pad}[2]{\frac{\partial{#1}}{\partial{#2}}}

\newcommand{\rpth}[1]{\left(#1\right)}
\newcommand{\spth}[1]{\left[#1\right]}
\newcommand{\cpth}[1]{\left\{#1\right\}}

\begin{document}

\title{Improving the scalability of a high-order atmospheric dynamics solver based on the \texttt{deal.II} library}

\author{Giuseppe Orlando$^{(1)}$, \\ 
        Tommaso Benacchio$^{(2)}$, Luca Bonaventura$^{(3)}$}

\date{}

\maketitle

\begin{center}
{
\small
$^{(1)}$  
CMAP, CNRS, \'{E}cole polytechnique, Institut Polytechnique de Paris \\ Route de Saclay, 91120 Palaiseau, France \\
{\tt giuseppe.orlando@polytechnique.edu} \\
\ \\
$^{(2)}$  
Weather Research, Danish Meteorological Institute \\
Sankt Kjelds Plads 11, 2100 Copenhagen, Denmark \\
{\tt tbo@dmi.dk} \\
\ \\
$^{(3)}$  
Dipartimento di Matematica, Politecnico di Milano \\
Piazza Leonardo da Vinci 32, 20133 Milano, Italy \\
{\tt luca.bonaventura@polimi.it} \\
}
\end{center}

\noindent

{\bf Keywords}: Numerical Weather Prediction, Non-conforming meshes, Flows over orography, IMEX schemes, \texttt{deal.II}

\pagebreak

\abstract{We present recent advances on the massively parallel performance of a numerical scheme for atmosphere dynamics applications based on the \texttt{deal.II} library. The implicit-explicit discontinuous finite element scheme is based on a matrix-free approach, meaning that no global sparse matrix is built and only the action of the linear operators on a vector is actually implemented. Following a profiling analysis, we focus on the performance optimization of the numerical method and describe the impact of different preconditioning and solving techniques in this framework. Moreover, we show how the use of the latest version of the \texttt{deal.II} library and of suitable execution flags can improve the parallel performance.}

\pagebreak

\section{Introduction}
\label{sec:intro}

Efficient numerical simulations of atmospheric flows are key to viable weather and climate forecasts and several related practical applications. Medium-range global numerical weather predictions (NWP) up to ten days ahead are typically required to complete within a one-hour operational threshold in the forecast cycles of meteorological centres. In addition, the computational meshes employed for these forecasts are being refined to head towards the km-scale, while experimental limited-area suites already reach the hectometric scale \cite{lean:2024}. Hence, in the current context of increasing demand of computational resources driven by increasingly high spatial resolutions, massively parallel model scalability is a fundamental requirement.

In this work, we present recent advances in the parallel performance of a Implicit-Explicit Discontinuous Galerkin (IMEX-DG) solver \cite{orlando:2022} employed for atmosphere dynamics applications \cite{orlando:2023, orlando:2024a, orlando:2024b}, following up on earlier analyses by the authors \cite{orlando:2025b}. The implementation is carried out in the framework of the \texttt{deal.II} library \cite{africa:2024, arndt:2023, bangerth:2007}. First, we show how the use of suitable preconditioning techniques can significantly improve the performance of the elliptic pressure solver. More specifically, we measure the impact of a novel preconditioning technique for the linear system associated with the pressure field that exploits the numerical formulation of the problem. Next, we assess the impact of activating execution flags and of using the most recent version of the library (9.6.2, \cite{africa:2024}). Finally, we present the result of employing a different strategy for the solution of a linear system.

The paper is structured as follows. In Section \ref{sec:model_num}, we briefly review the model equations and the relevant details of the numerical methodology used in this work. The profiling and parallel performance analysis in numerical experiments with a three-dimensional benchmark of atmospheric flow are presented in Section \ref{sec:test}. Finally, some conclusions are reported in Section \ref{sec:conclu}.

\section{The model equations and some details of the numerical method}
\label{sec:model_num}

The compressible Euler equations of gas dynamics represent the most comprehensive mathematical model for atmospheric flows \cite{steppeler:2003}. Considering for simplicity the dry, non-rotating case, the mathematical model reads as follows:
\begin{eqnarray}\label{eq:euler_comp}
   &&\pad{\rho}{t} + \dive\rpth{\rho\vel} = 0 \nonumber \\
   &&\pad{\rpth{\rho\vel}}{t} + \dive\rpth{\rho\vel \otimes \vel} + \grad p = \rho\mathbf{g} \\
   &&\pad{\rpth{\rho E}}{t} + \dive\spth{\rpth{\rho h + \rho k}\vel} = \rho\mathbf{g} \cdot \vel, \nonumber
\end{eqnarray}
where $t \in (0, T_{f}]$ is the temporal coordinate, supplied with suitable initial and boundary conditions on the computational domain. Here $\otimes$ denotes the tensor product, $T_{f}$ is the final time, $\rho$ is the density, $\mathbf{u}$ is the fluid velocity, and $p$ is the pressure. Moreover, $\mathbf{g} = -g\mathbf{k}$ is the acceleration of gravity, with $g = \SI{9.81}{\meter\per\second\squared}$ and $\mathbf{k}$ being the upward pointing unit vector in the standard Cartesian frame of reference. Finally, $E$ denotes the total energy per unit of mass, $h$ is the specific enthalpy and $k = 1/2\left|\vel\right|^{2}$ is the specific kinetic energy. System \eqref{eq:euler_comp} is complemented by the equation of state of ideal gases, given by $\rho e = \frac{1}{\gamma - 1}p$, where $e$ denotes the specific internal energy $e = E - k$ and the isentropic exponent $\gamma$ is taken equal to $1.4$.

The time discretization is based on an Implicit-Explicit (IMEX) Runge-Kutta (RK) method \cite{kennedy:2003}, which is a widely employed approach for ODE systems that include both stiff and non-stiff components. IMEX-RK schemes are represented compactly by the companion Butcher tableaux:
\begin{center}
    \begin{tabular}{c|c}
	$\mathbf{c}$ & $\mathbf{A}$ \\
	\hline
	& $\mathbf{b}^{\top}$
    \end{tabular}
    \qquad
    \begin{tabular}{c|c}
	$\tilde{\mathbf{c}}$ & $\tilde{\mathbf{A}}$  \\
	\hline
	& $\tilde{\mathbf{b}}^{\top}$
    \end{tabular}
\end{center}
with $\mathbf{A} = \cpth{a_{lm}}, \mathbf{b} = \cpth{b_{l}}, \mathbf{c} = \cpth{c_{l}}$, $l,m = 1\dots s$, denoting the coefficients of the explicit method, $\tilde{\mathbf{A}} = \cpth{\tilde{a}_{lm}}, \tilde{\mathbf{b}} = \cpth{\tilde{b}_{l}}$, and $\tilde{\mathbf{c}} = \cpth{\tilde{c}_{l}}$, $l,m = 1\dots s$ denoting the coefficients of the implicit method, and $s$ representing the number of stages of the method. We refer, e.g., to \cite{kennedy:2003, pareschi:2005} for a detailed analysis of the order and coupling conditions of the two companion schemes. A generic IMEX-RK $l$ stage with time step $\Delta t$ for the Euler equations reads therefore as follows:
\begin{eqnarray}\label{eq:stage_euler}
   &&\rho^{(l)} = \rho^{n} - \sum_{m=1}^{l-1}a_{lm}\Delta t\dive\rpth{\rho^{(m)}\vel^{(m)}} \nonumber \\
   &&\rho^{(l)}\mathbf{u}^{(l)} + \tilde{a}_{ll}\Delta t\grad p^{(l)} = \mathbf{m}^{(l)} \\
   &&\rho^{(l)}E^{(l)} + \tilde{a}_{ll}\Delta t\dive\rpth{h^{(l)}\rho^{(l)}\vel^{(l)}} = \hat{e}^{(l)}, \nonumber
\end{eqnarray}
where we have set
\begin{eqnarray}
    \mathbf{m}^{(l)} &=& \rho^{n}\vel^{n} - \sum_{m=1}^{l-1}a_{lm}\Delta t\dive\rpth{\rho^{(m)}\vel^{(m)} \otimes \vel^{(m)}} - \sum_{m=1}^{l-1}\tilde{a}_{lm}\Delta t\grad p^{(m)} \\
    \hat{\mathbf{e}}^{(l)} &=& \rho^{n}E^{n} - \sum_{m=1}^{l-1}\tilde{a}_{lm}\Delta t\dive\rpth{h^{(m)}\rho^{(m)}\vel^{(m)}}) - \sum_{m=1}^{l-1}a_{lm}\Delta t \dive\rpth{k^{(m)}\rho^{(m)}\vel^{(m)}}
\end{eqnarray}
Notice that, substituting formally $\rho^{(l)}\vel^{(l)}$ into the energy equation and taking into account the definitions $\rho E = \rho e + \rho k$ and $h = e + p/\rho$, the following nonlinear Helmholtz-type equation for the pressure is obtained:
\begin{eqnarray}\label{eq:Helmholtz}
    \rho^{(l)}\spth{e(p^{(l)},\rho^{(l)}) + k^{(l)}} &-& \tilde{a}_{ll}^{2}\Delta t^{2}\dive\spth{\rpth{e(p^{(l)},\rho^{(l)}) + \frac{p^{(l)}}{\rho^{(l)}}}\nabla p^{(l)}} \nonumber \\
    &+& \tilde{a}_{ll}\Delta t\dive\spth{\rpth{e(p^{(l)},\rho^{(l)}) + \frac{p^{(l)}}{\rho^{(l)}}}\mathbf{m}^{(l)}}
    = \hat{e}^{(l)}, 
\end{eqnarray}
which is solved through a fixed point procedure. We refer to \cite{orlando:2022, orlando:2025a, orlando:2025c} for further references and details on the theoretical properties of the method.

The spatial discretization is based on a Discontinuous Galerkin (DG) method \cite{giraldo:2020}, which combines high-order accuracy and flexibility in a highly data-local framework. In particular, this method is particularly well suited for mesh adaptive approaches \cite{orlando:2023, orlando:2024a}. The shape functions correspond to the products of Lagrange interpolation polynomials for the support points of $\left(r + 1\right)$-order Gauss-Lobatto quadrature rule in each coordinate direction, with $r$ denoting the polynomial degree. The discrete formulation of method \eqref{eq:stage_euler} can therefore be expressed as
\begin{eqnarray}
   &&\mathbf{A}^{(l)}\mathbf{U}^{(l)} + \mathbf{B}^{(l)}\mathbf{P}^{(l)} = \mathbf{F}^{(l)} \label{eq:discr_form_momentum} \\
   &&\mathbf{C}^{(l)}\mathbf{U}^{(l)} + \mathbf{D}^{(l)}\mathbf{P}^{(l)} = \mathbf{G}^{(l)},
\end{eqnarray}
where $\mathbf{U}^{(l)}$ represents the vector of degrees of freedom of the velocity field, while $\mathbf{P}^{(l)}$ is the vector of the degrees of freedom of the pressure field. We refer to \cite{orlando:2023, orlando:2025a, orlando:2025b} for the detailed expression of all the matrices and vectors. We report only the expression of the matrix $A^{(l)}$ that will be useful for the discussion in Section \ref{sec:test}:
\begin{equation}\label{eq:modified_mass_matrix_velocity}
    A_{ij}^{(l)} = \sum_{K \in \mathcal{T}_{\mathcal{H}}} \int_{K} \rho^{(l)}\boldsymbol{\varphi}_{j} \cdot \boldsymbol{\varphi}_{i}\mathrm{d}\mathbf{x},
\end{equation}
where $\mathcal{T}_{\mathcal{H}}$ represents the computational mesh, $\boldsymbol{\varphi}_{i}$ denotes the basis function of the space of polynomial functions employed to discretize the velocity, and $\mathrm{d}\mathbf{x}$ is the spatial element of integration. Formally, from \eqref{eq:discr_form_momentum} one can derive
\begin{equation}
    \mathbf{U}^{(l)} = \left(\mathbf{A}^{(l)}\right)^{-1}\left(\mathbf{F}^{(l)} - \mathbf{B}^{(l)}\mathbf{P}^{(l)}\right), 
\end{equation}
and obtain
\begin{equation}\label{eq:fixed_point_discrete}
    \mathbf{D}^{(l)}\mathbf{P}^{(l)} + \mathbf{C}^{(l)}\left(\mathbf{A}^{(l)}\right)^{-1}\left(\mathbf{F}^{(l)} - \mathbf{B}^{(l)}\mathbf{P}^{(l)}\right) = \mathbf{G}^{(l)}.
\end{equation}
The above system is then solved following the fixed point procedure described in \cite{orlando:2022}. In Sections \ref{ssec:preconditinoer} and \ref{ssec:block-diag} we describe the impact of using different preconditioning and solving techniques for system \eqref{eq:fixed_point_discrete}.

\section{Numerical results}
\label{sec:test}

In this Section, we show results of simulations of an idealized three-dimensional test case of atmospheric flow over orography \cite{melvin:2019, orlando:2024a} using the numerical method described in the previous Section and focusing on its scalability. The simulations have been run using up to 1024 2x AMD EPYC Rome 7H12 64c\@ 2.6GHz CPUs at MeluXina HPC facility\footnote{\url{https://docs.lxp.lu/}} and OpenMPI 4.1.5 has been employed. The compiler is GCC version 12.3 and the Vectorization level is 256 bits.

We consider a three-dimensional configuration of a flow over a bell-shaped hill, already studied in \cite{melvin:2019, orlando:2023, orlando:2024a}. The computational domain is $\SI[parse-numbers=false]{\left(0, 60\right) \times \left(0, 40\right) \times \left(0, 16\right)}{\kilo\meter}$. The mountain profile is defined as follows:
\begin{equation}
    h(x,y) =  h_{c}{\left[1 + \left(\frac{x - x_{c}}{a_{c}}\right)^{2} + \left(\frac{y - y_{c}}{a_{c}}\right)^{2}\right]^{-\frac{3}{2}}},
\end{equation}
with $h_{c} = \SI{400}{\meter}, a_{c} = \SI{1}{\kilo\meter}, x_{c} = \SI{30}{\kilo\meter},$ and $y_{c} = \SI{20}{\kilo\meter}$. The buoyancy frequency is $N = \SI{0.01}{\per\second}$, whereas the background velocity is $\overline{u} = \SI{10}{\meter\per\second}$. The final time is $T_{f} = \SI{1}{\hour}$. The initial conditions read as follows \cite{benacchio:2014}:
\begin{equation}
    p = p_{ref}\left\{1 - \frac{g}{N^{2}}\Gamma\frac{\rho_{ref}}{p_{ref}}\left[1 - \exp\left(-\frac{N^{2}z}{g}\right)\right]\right\}^{1/\Gamma}, \quad \rho = \rho_{ref}\left(\frac{p}{p_{ref}}\right)^{1/\gamma}\exp\left(-\frac{N^{2}z}{g}\right),
\end{equation}
where $p_{ref} = \SI[parse-numbers = false]{10^{5}}{\pascal}$ and $\rho_{ref} = \frac{p_{ref}}{R T_{ref}}$, with $T_{ref} = \SI{293.15}{\kelvin}$. Finally, we set $\Gamma = \frac{\gamma - 1}{\gamma}$. We refer to \cite{orlando:2024a} for the implementation of boundary conditions. Unless differently stated, we consider a) a uniform mesh composed of $N_{el} = 120 \times 80 \times 32 = 307200$ elements with polynomial degree $r = 4$, leading to about $38.5$ million unknowns for each scalar variable and a resolution of $\SI{125}{\meter}$;  and b) a non-conforming mesh composed of $N_{el} = 204816$ elements, yielding around $25.6$ millions of unknowns for each scalar variable and a maximum resolution of $\SI{62.5}{\meter}$. A non-conforming mesh is characterized by neighbouring cells with different resolution on both the horizontal and the vertical direction \cite{orlando:2024a, orlando:2025b} and its use provides sizeable computational savings \cite{orlando:2024a, orlando:2025b}. The following optimization flags are employed throughout the runs:
\\~\\
\texttt{-O2 -funroll-loops -funroll-all-loops -fstrict-aliasing -Wno-unused-local-typedefs}
\\~\\
which are the standard ones in the Release mode of \texttt{deal.II} (see \url{https://www.dealii.org/developer/users/cmake_user.html}). Unfortunately, due to limitations on computational resources, the scaling analyses reported in this Section could only be performed once for each of the different configurations.

\subsection{Impact of preconditioning technique}
\label{ssec:preconditinoer}

The preconditioned conjugate gradient method with a geometric multigrid preconditioner is employed to solve the symmetric positive definite linear systems, while the GMRES solver with a Jacobi preconditioner is employed for the solution of non-symmetric linear systems. A matrix-free approach is employed \cite{arndt:2023}, meaning that no global sparse matrix is built and only the action of the linear operators on a vector is actually implemented. This strategy is very efficient for high-order discretization methods \cite{kronbichler:2019}. However, it leads to some difficulties in the preconditioning, since the matrix is not available and standard techniques as ILU cannot be employed. For symmetric positive-definite systems, an efficient preconditioning technique compatible with the matrix-free approach is based on the so-called Chebyshev polynomial that relies on an estimate of the eigenvalues of the matrix so as to damp the eigenvalue range \cite{lanczos:1950}.

The approach based on the eigenvalue estimate cannot be directly employed for the solution of \eqref{eq:fixed_point_discrete} because it requires symmetry and positive definiteness of the preconditioned matrix \cite{varga:2009}. Moreover, the entries of the matrices $(\mathbf{A}^{(l)})^{-1}$, $\mathbf{D}^{(l)}$, and $\mathbf{C}^{(l)}$ are not readily available because of the matrix-free approach. In a matrix-free framework, a matrix representation of the operator can be obtained by applying the operator on all unit vectors. Apparently, this is a very inefficient procedure because it requires to perform $n$ operator evaluations for a $n \times n$ matrix. In practice, the integration is so efficient that the computation completes without significant overhead \footnote{\url{https://www.dealii.org/current/doxygen/deal.II/step_37.html}}. However, using this procedure for the computation of $(\mathbf{A}^{(l)})^{-1}$, $\mathbf{D}^{(l)}$, and $\mathbf{C}^{(l)}$ turns out to be very inefficient because it requires the evaluation of three operators as well as the solution of a linear system at each operator evaluation to compute $(\mathbf{A}^{(l)})^{-1}$. Hence, the simplest solution consists in considering only the contribution provided by $\mathbf{D}^{(l)}$, which takes into account the internal energy. For low Mach number flows, as those considered in this work, this approximation can be considered reliable since velocities are typically low and the internal energy is therefore dominant with respect to the kinetic energy. This is the strategy adopted, e.g., in \cite{orlando:2024a}, and we refer to it as \textit{internal energy preconditioner}.

We consider here a novel alternative approach that can be easily implemented in a matrix-free framework and consists of defining an operator that takes into account the underlying nonlinear Helhmholtz-type equation \eqref{eq:Helmholtz} when solving \eqref{eq:fixed_point_discrete}. More specifically, we consider the volumetric contribution of the weak form associated to a DG discretization of \eqref{eq:Helmholtz}. Hence, the preconditioner is the inverse of the diagonal part of the matrix $\tilde{\mathbf{S}}^{(l)}$ whose entries are
\begin{equation}\label{eq:entries_Helholtz_precon}
    \tilde{S}^{(l)}_{ij} = \sum_{K \in \mathcal{T}_{\mathcal{H}}}\int_{K}\rho^{(n,l)}e\rpth{\Psi_{j}, \rho^{(l)}}\Psi_{i}\mathrm{d}\mathbf{x} + \sum_{K \in \mathcal{T}_{\mathcal{H}}}\int_{K}\tilde{a}_{ll}^{2}\Delta t^{2}\rpth{e\rpth{p^{(l)},\rho^{(l)}} + \frac{p^{(l)}}{\rho^{(l)}}}\grad\Psi_{j} \cdot \grad\Psi_{i}\mathrm{d}\mathbf{x}.
\end{equation}
where $\Psi_{i}$ denotes the basis function of the space of polynomial functions employed to discretize the pressure. Matrix $\tilde{\mathbf{S}}^{(l)}$ is the sum of two contributions: a first one that considers the internal energy and that coincides with matrix $\mathbf{D}^{(l)}$, and a second one that considers the contribution due to specific enthalpy and takes into account the underlying elliptic character of \eqref{eq:Helmholtz}. The pressure $p^{(l)}$ in \eqref{eq:entries_Helholtz_precon} in the fixed point procedure is evaluated at the previous fixed point iteration. We refer to this strategy as \textit{Helmholtz preconditioner}.

\begin{figure}[h!]
    \centering
    \begin{subfigure}{0.475\textwidth}
	\centering
        \includegraphics[width = 0.95\textwidth]{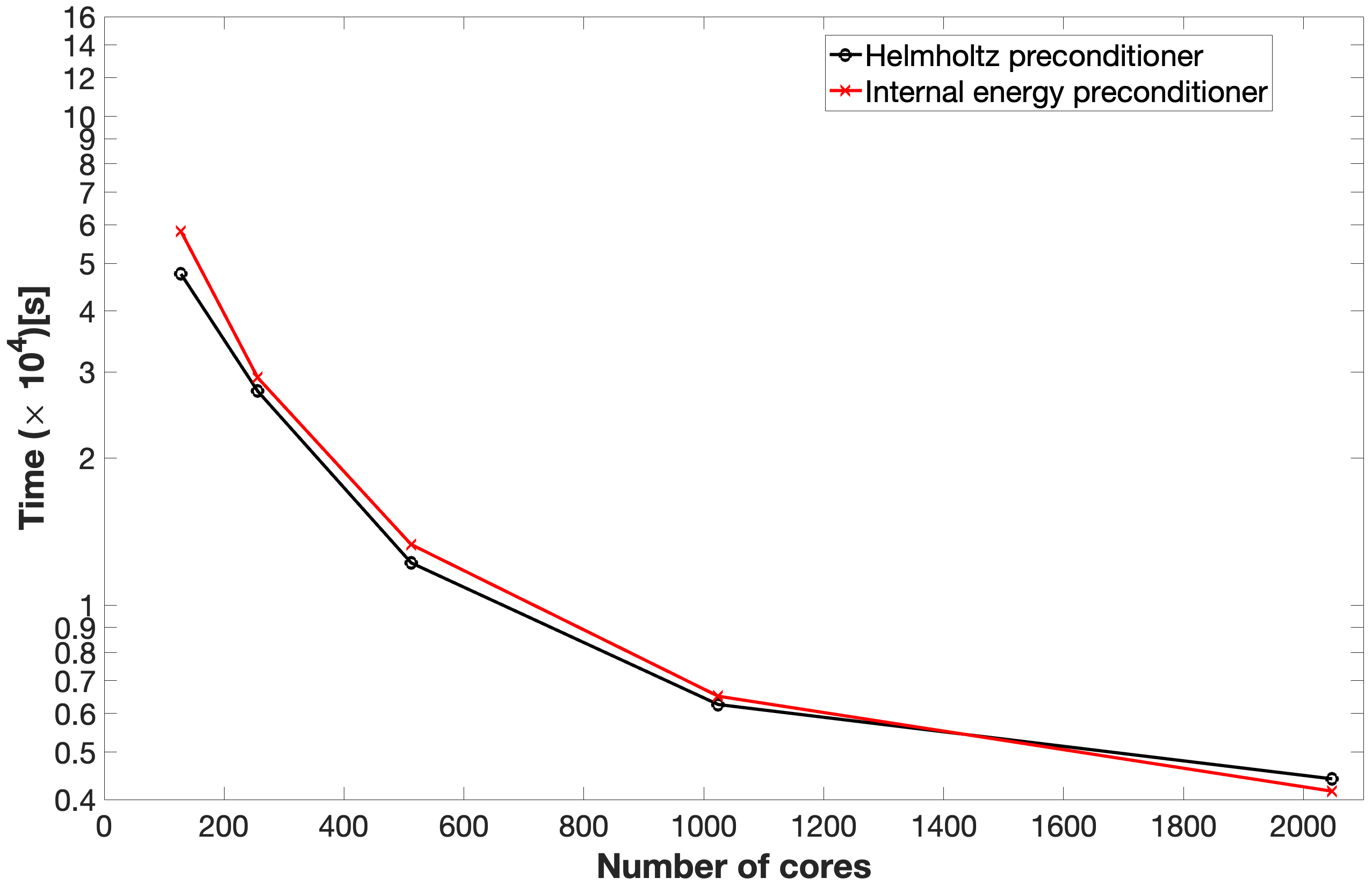}
    \end{subfigure}
    \begin{subfigure}{0.475\textwidth}
	\centering
        \includegraphics[width = 0.95\textwidth]{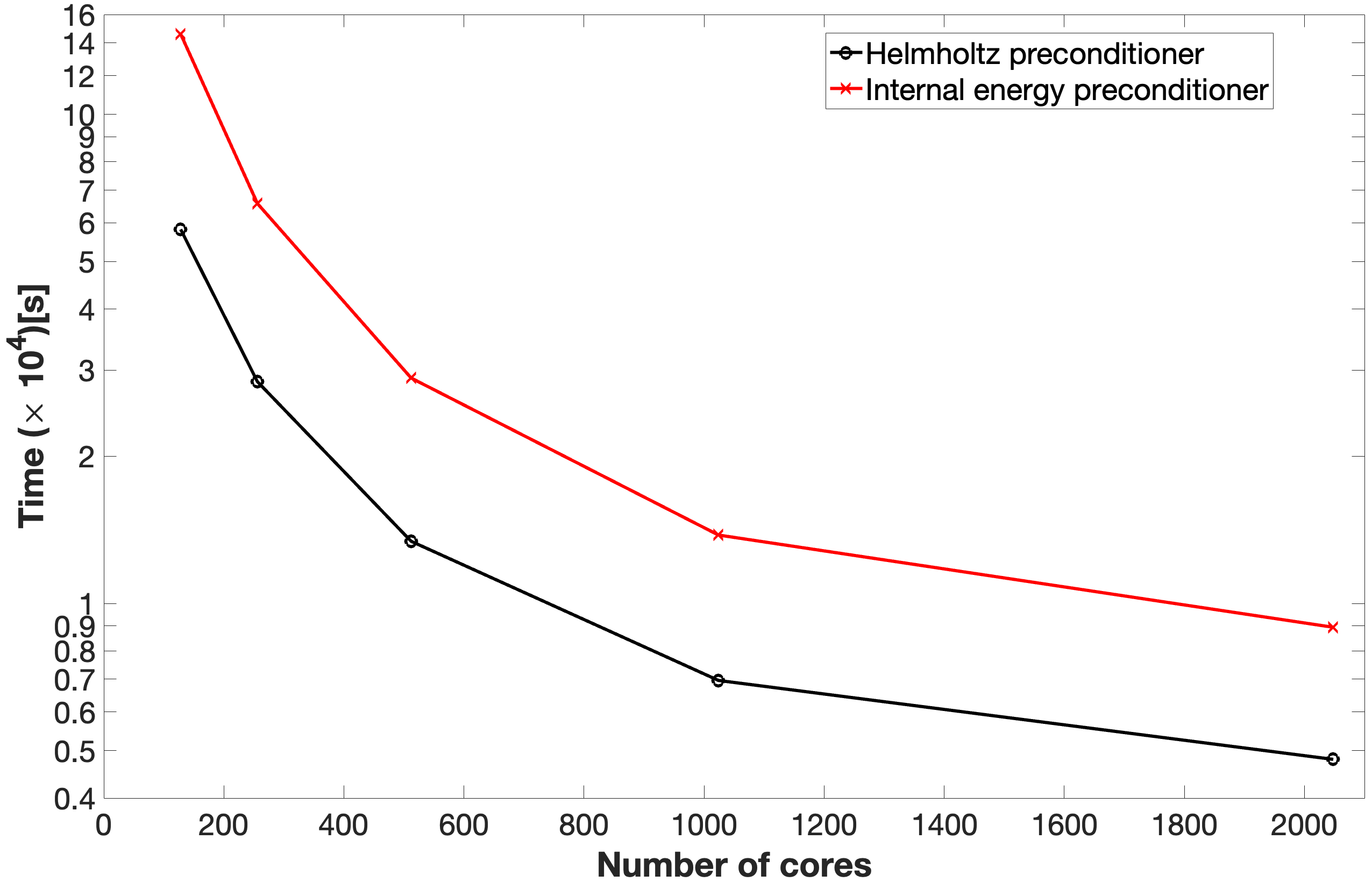}
    \end{subfigure}
    \caption{Analysis of preconditioning technique, wall-clock time in the numerical solution of problem \eqref{eq:stage_euler}, flow over a 3D hill test case. Left: uniform mesh. Right: non-conforming mesh. The black lines show the results with the \textit{Helmholtz preconditioner} approach, whereas the red lines show the results with \textit{internal energy preconditioner} approach.}
    \label{fig:WT_precon}
\end{figure}

An analysis of time to solution using the two preconditioning techniques reveals that, when using a uniform mesh, the use of the Helmholtz preconditioner gives a relatively modest computational time saving - about $18\%$ with 128 cores (1 node) - that decreases for higher core counts (Figure \ref{fig:WT_precon}). A more sizeable computational time saving is established for the non-conforming mesh which amounts to $60\%$ employing 128 cores and to $46 \%$ employing 2048 cores (Figure \ref{fig:WT_precon}). Hence, strategy \textit{Helmholtz preconditioner} provides a significant advantage in terms of time-to-solution, in particular for non-conforming meshes. In the following, we will restrict our attention only to this preconditioning technique and to experiments using the non-conforming mesh.

\subsection{Impact of execution flags and version}
\label{ssec:flags}

\begin{figure}[h!]
    \centering
    \begin{subfigure}{0.475\textwidth}
	\centering
        \includegraphics[width = 0.95\textwidth]{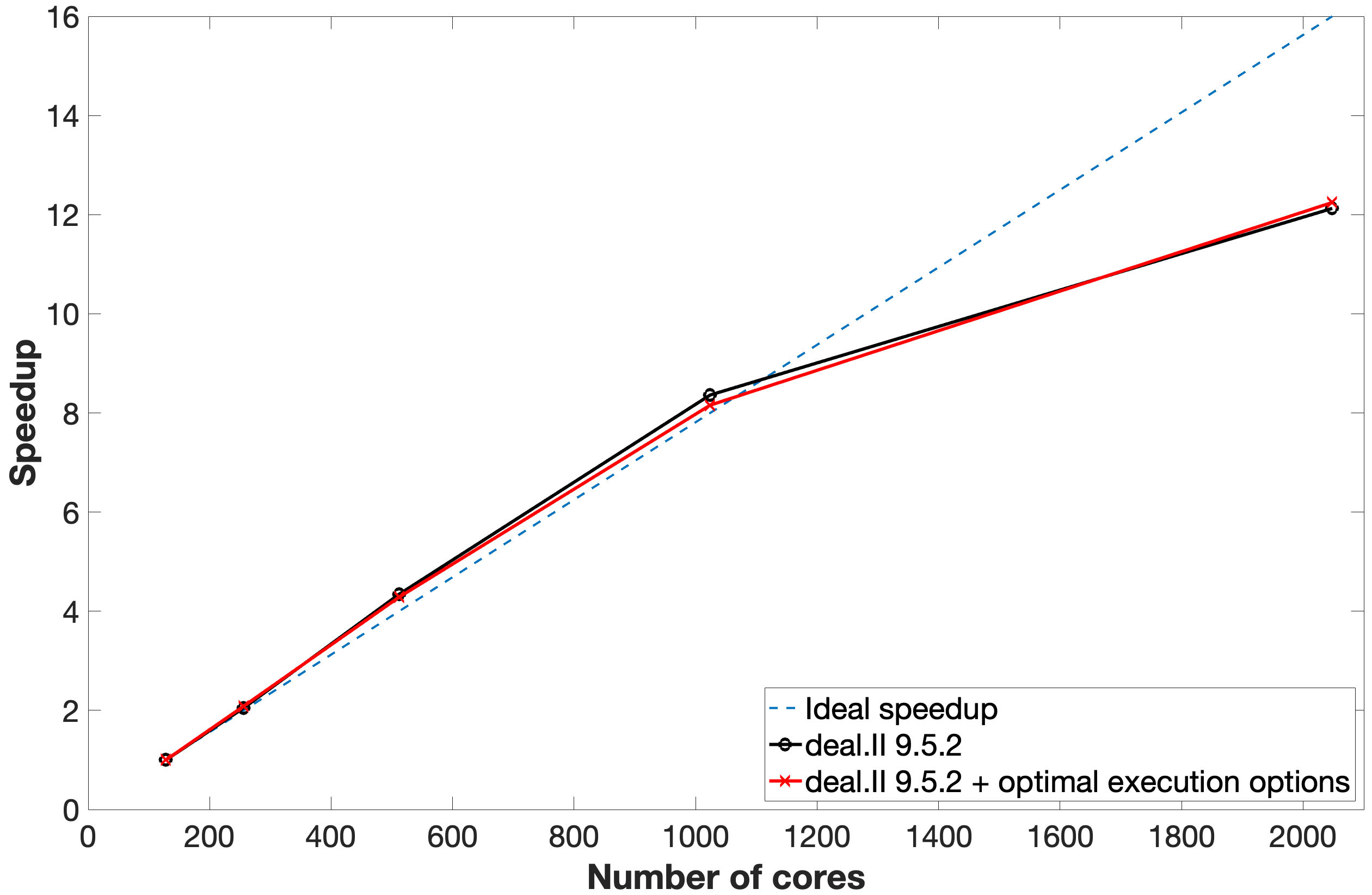}
    \end{subfigure}
    \begin{subfigure}{0.475\textwidth}
	\centering
        \includegraphics[width = 0.95\textwidth]{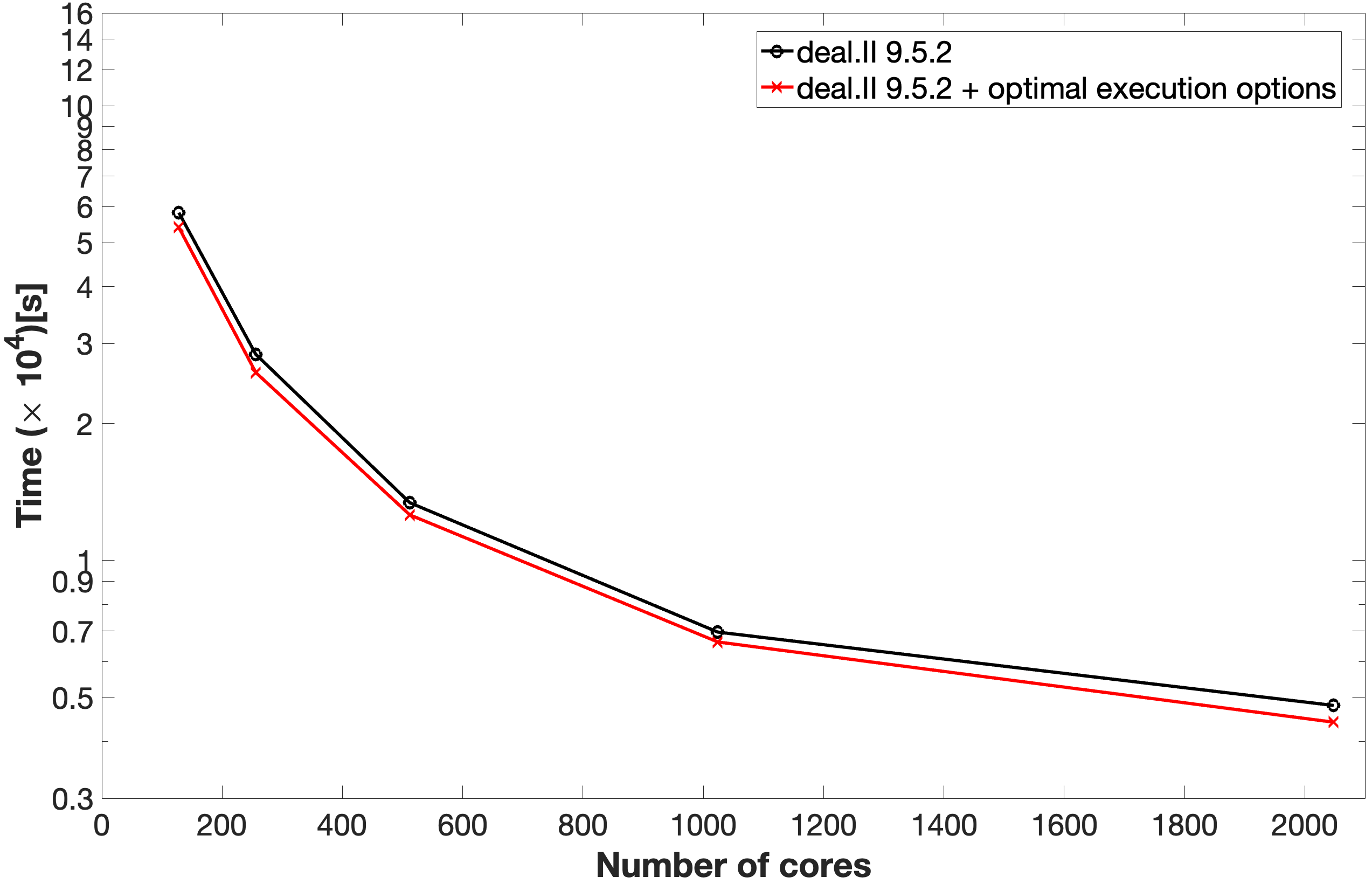}
    \end{subfigure}
    \caption{Performance impact of using different execution flags in the numerical solution of problem \eqref{eq:stage_euler}, flow over a 3D hill test case. Left: strong scaling. Right: WT time. The results are obtained with the non-conforming mesh. The black lines show the results with the standard execution flags, whereas the red lines show the results with optimal execution flags.}
    \label{fig:flags}
\end{figure}

Next, we analyze the impact of some execution flags and of the use of the \texttt{deal.II} 9.6.2 version \cite{africa:2024} instead of the \texttt{deal.II} 9.5.2 version \cite{arndt:2023}. Thanks to the support of EPICURE team at MeluXina, it has been shown that the execution flags
\\~\\
\texttt{export DEAL\_II\_NUM\_THREADS="\$SLURM\_CPUS\_PER\_TASK"} \\
\texttt{export OMP\_NUM\_THREADS="\$SLURM\_CPUS\_PER\_TASK"}
\\~\\
improve the parallel performance of the solver in terms of both strong scaling and wall-clock time (Figure \ref{fig:flags}). However, an even more sizeable improvement in the parallel performance is obtained using the \texttt{deal.II} 9.6.2 version, leading to a computational time saving of around $28\%$ compared to using the 9.5.2 version in a run with 2048 cores (i.e. 16 nodes). For the results in the following Section \ref{ssec:block-diag}, \texttt{deal.II} 9.6.2 version with the optimal execution flags is therefore employed.

\begin{figure}[h!]
    \centering
    \begin{subfigure}{0.475\textwidth}
	\centering
        \includegraphics[width = 0.95\textwidth]{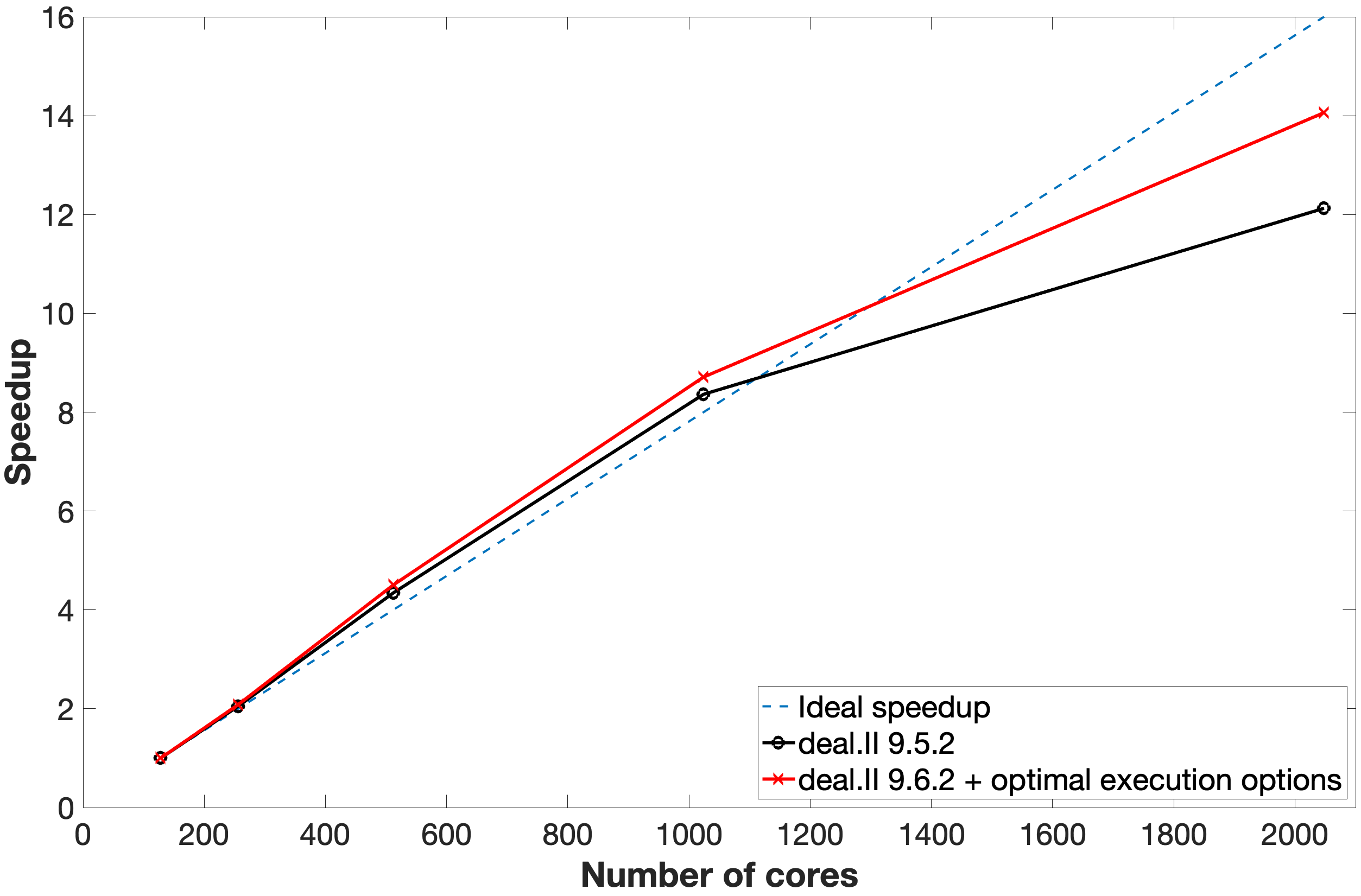}
    \end{subfigure}
    \begin{subfigure}{0.475\textwidth}
	\centering
        \includegraphics[width = 0.95\textwidth]{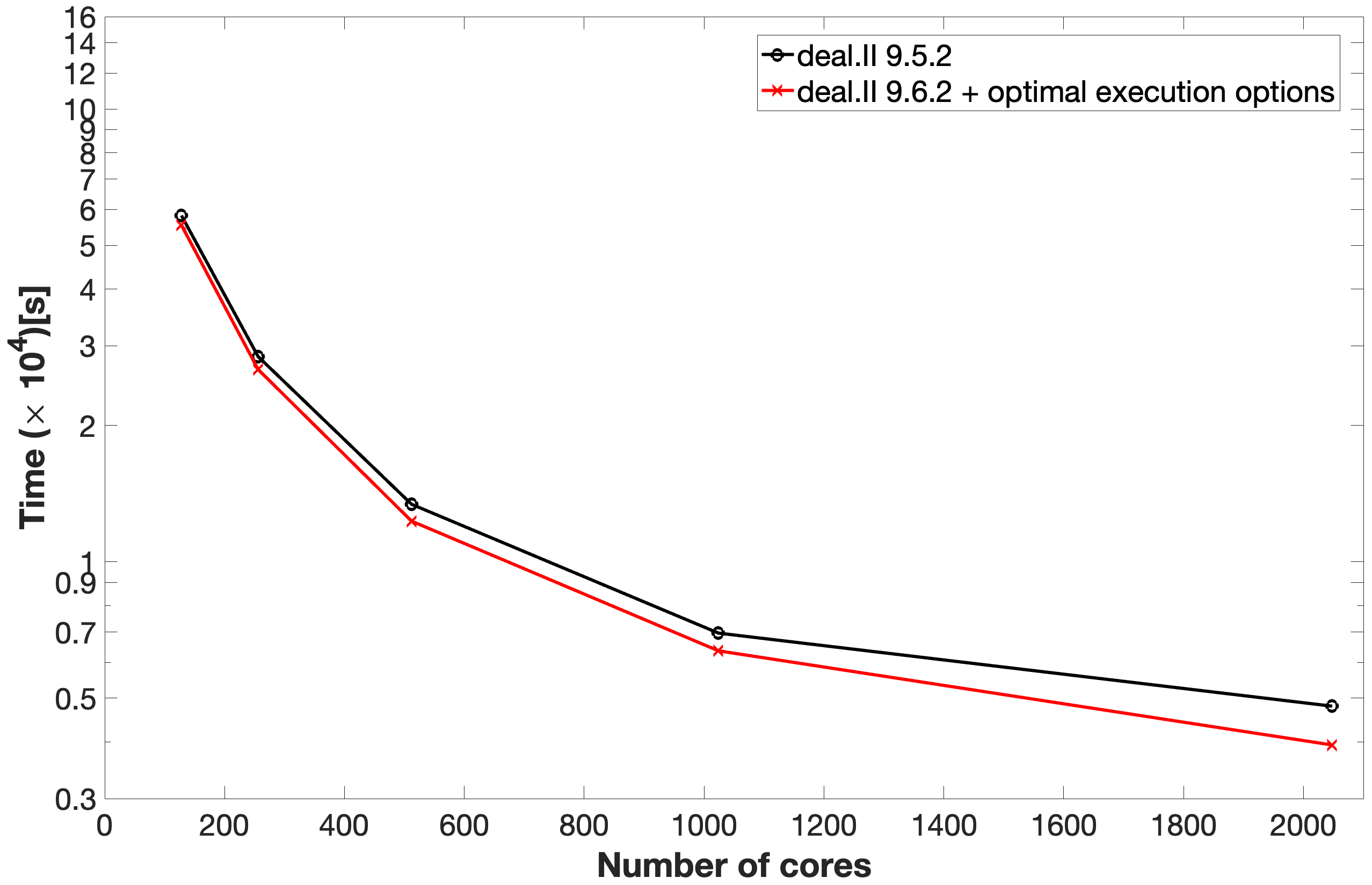}
    \end{subfigure}
    \caption{Performance impact of using different \texttt{deal.II} versions in the numerical solution of problem \eqref{eq:stage_euler}, flow over a 3D hill test case. Left: strong scaling. Right: WT time. The results are obtained with the non-conforming mesh. The black lines show the results obtained using the \texttt{deal.II} 9.6.2 version and with optimal execution flags, whereas the red lines show the results obtained using the \texttt{deal.II} 9.5.2 version with standard execution flags.}
    \label{fig:version_uniform}
\end{figure}

\subsection{Fast evaluation of block-diagonal operators}
\label{ssec:block-diag}

\begin{figure}[h!]
    \centering
    \includegraphics[width = 0.7\textwidth]{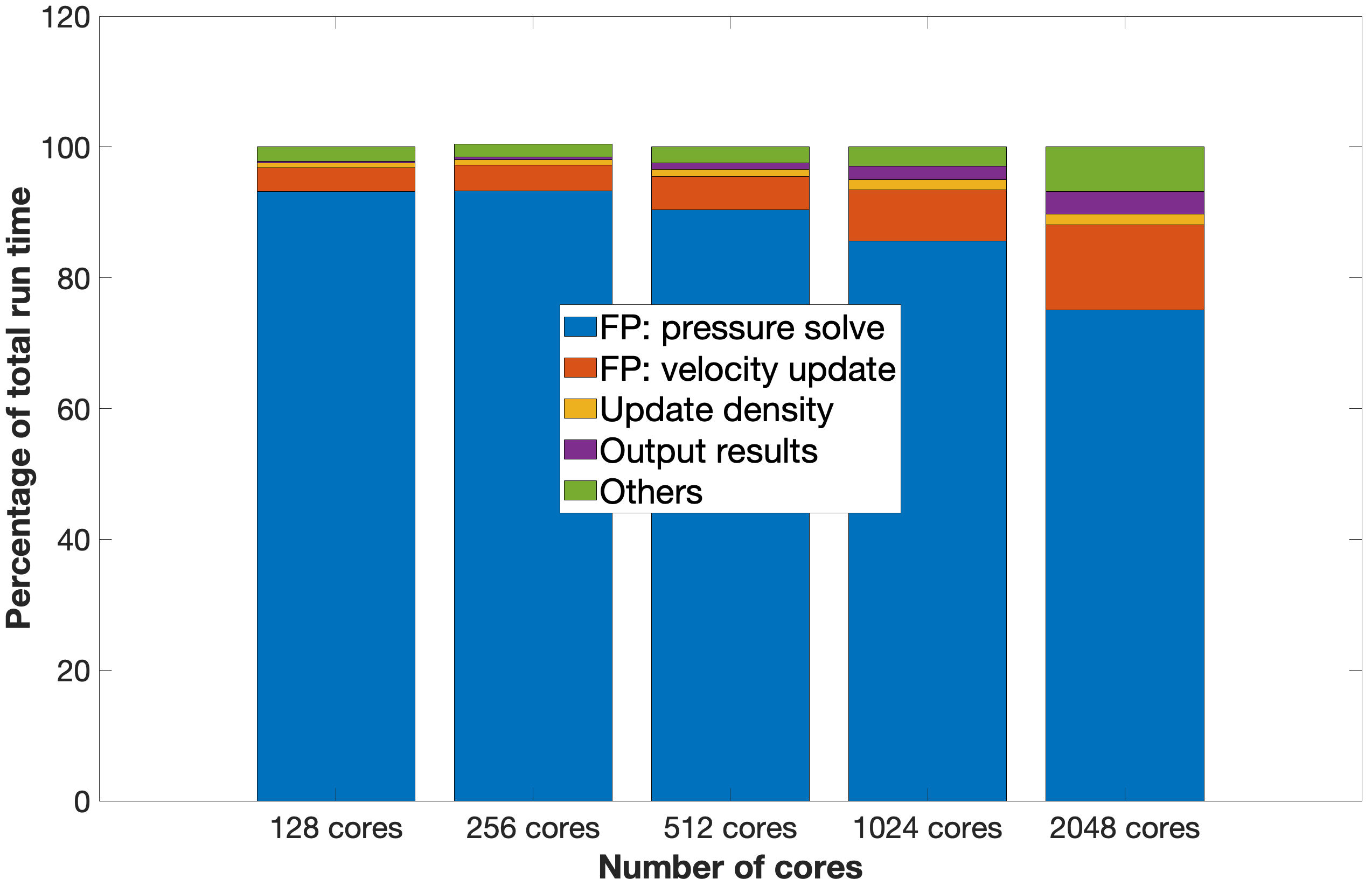}
    \caption{Distribution of the computational time spent in various blocks of the algorithm for the solution of problem \eqref{eq:stage_euler}, flow over a 3D hill test case as a function of number of cores.}
    \label{fig:profiling_consistent_integration}
\end{figure}

Finally, we have investigated the distribution of the computational time spent in the different part of the algorithms. A profiling analysis show that, as expected, the vast majority of execution time is spent in solving problem \eqref{eq:fixed_point_discrete} to compute the pressure field (Figure \ref{fig:profiling_consistent_integration}). The results are obtained using the non-conforming mesh, but analogous considerations hold for the uniform mesh (not shown). Apart from the computational time spent in computing the pressure field, one can immediately notice an increasing relative cost of computing the velocity field (orange bar in Figure \ref{fig:profiling_consistent_integration}) which ranges from about $3.6\%$ with 128 cores to about $13\%$ with 2048 cores. This part of the algorithm consists in computing $(\mathbf{A}^{(l)})^{-1}$ (we refer to \cite{orlando:2022, orlando:2023} for all the details), hence solving a linear system. However, this operation is repeated several times in the fixed point loop because of the matrix-free approach, which explains the increasing computational burden. 

The matrix $\mathbf{A}^{(l)}$ \eqref{eq:modified_mass_matrix_velocity} is a simple block-diagonal operator. Efficient techniques to evaluate and to invert block-diagonal operators in a DG framework were presented in \cite{kronbichler:2019}. However, they require the use of the same number of quadrature points and degrees of freedom along each coordinate direction, since this yields, after a change of basis, diagonal blocks in the algebraic structure of the the block-diagonal operator \cite{kronbichler:2019}. More specifically, as explained in \cite{kronbichler:2019, kronbichler:2016}, the block-diagonal operator is inverted as $\mathbf{S}^{-T} \mathbf{J}^{-1} \mathbf{S}^{-1}.$ Here $\mathbf{J}$ is a diagonal matrix, whose entries are equal to the determinant of the Jacobian times the quadrature weight, while $\mathbf{S}$ is a square matrix with basis functions in the row and quadrature points in columns, i.e. $S_{ij} = \varphi_{i}(\mathbf{x}_{j})$. The matrix $\mathbf{S}$ is then constructed as the Kronecker product (tensor product) of small one-dimensional matrices that can be inverted efficiently at each stage. In a nodal DG framework with nodes located at $r+1$ on the quadrature points of a Gauss-Legendre-Lobatto formula, as the one considered here, this means that it is possible to employ a Gauss-Legendre formula with $r+1$ quadrature points along each coordinate direction for the numerical integration.

The main drawback of this approach is that an aliasing error is introduced. Since numerical integration based on the Gauss-Legendre formula with $r+1$ quadrature points integrates exactly polynomials up to degree $2r + 1$, the strategy presented in \cite{kronbichler:2019} does not introduce aliasing errors for classical mass matrices that consist of the product of two basis function. However, the matrix $\mathbf{A}^{(l)}$ \eqref{eq:modified_mass_matrix_velocity} is a modified mass matrix and consists of the product between three polynomials. Hence, an aliasing error is introduced for $r > 1$. The numerical integration for the results shown so far is based on the so-called consistent integration or over-integration (see \cite{orlando:2024b} for further details). More specifically, a Gauss-Legendre formula with $2r+1$ quadrature points along each coordinate direction is employed. Notice that the choice to consider a polynomial representation for density, velocity, and pressure avoids any polynomial division in the case of the ideal gas law and therefore all the integrals can be computed without aliasing errors thanks to the aforementioned consistent integration. Hence, the choice of the numerical quadrature formula is related to the need to integrate all the terms so as to avoid aliasing errors. 

\begin{figure}[h!]
    \centering
    \begin{subfigure}{0.475\textwidth}
	\centering
        \includegraphics[width = 0.95\textwidth]{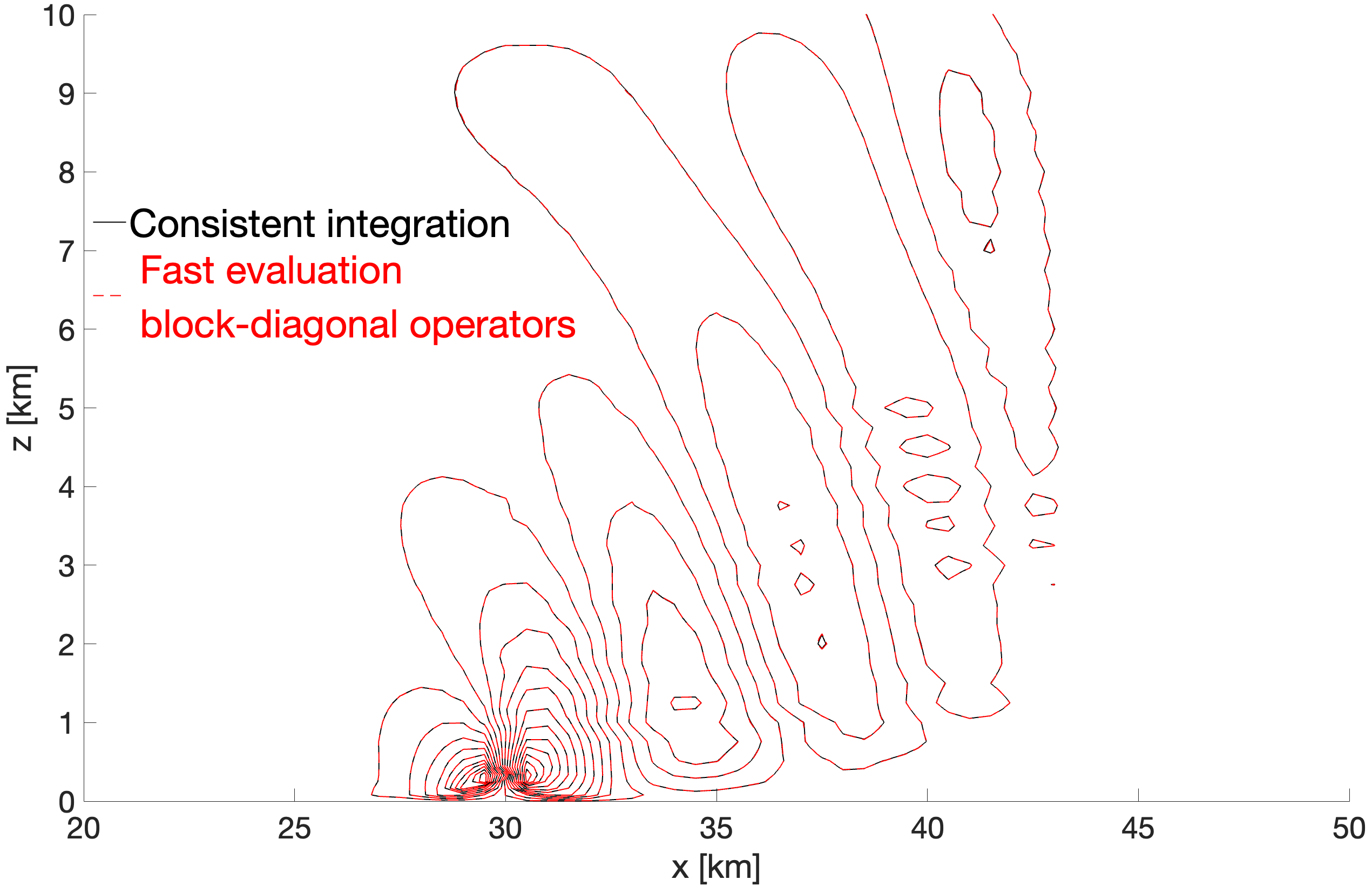}
    \end{subfigure}
    \begin{subfigure}{0.475\textwidth}
	\centering
        \includegraphics[width = 0.95\textwidth]{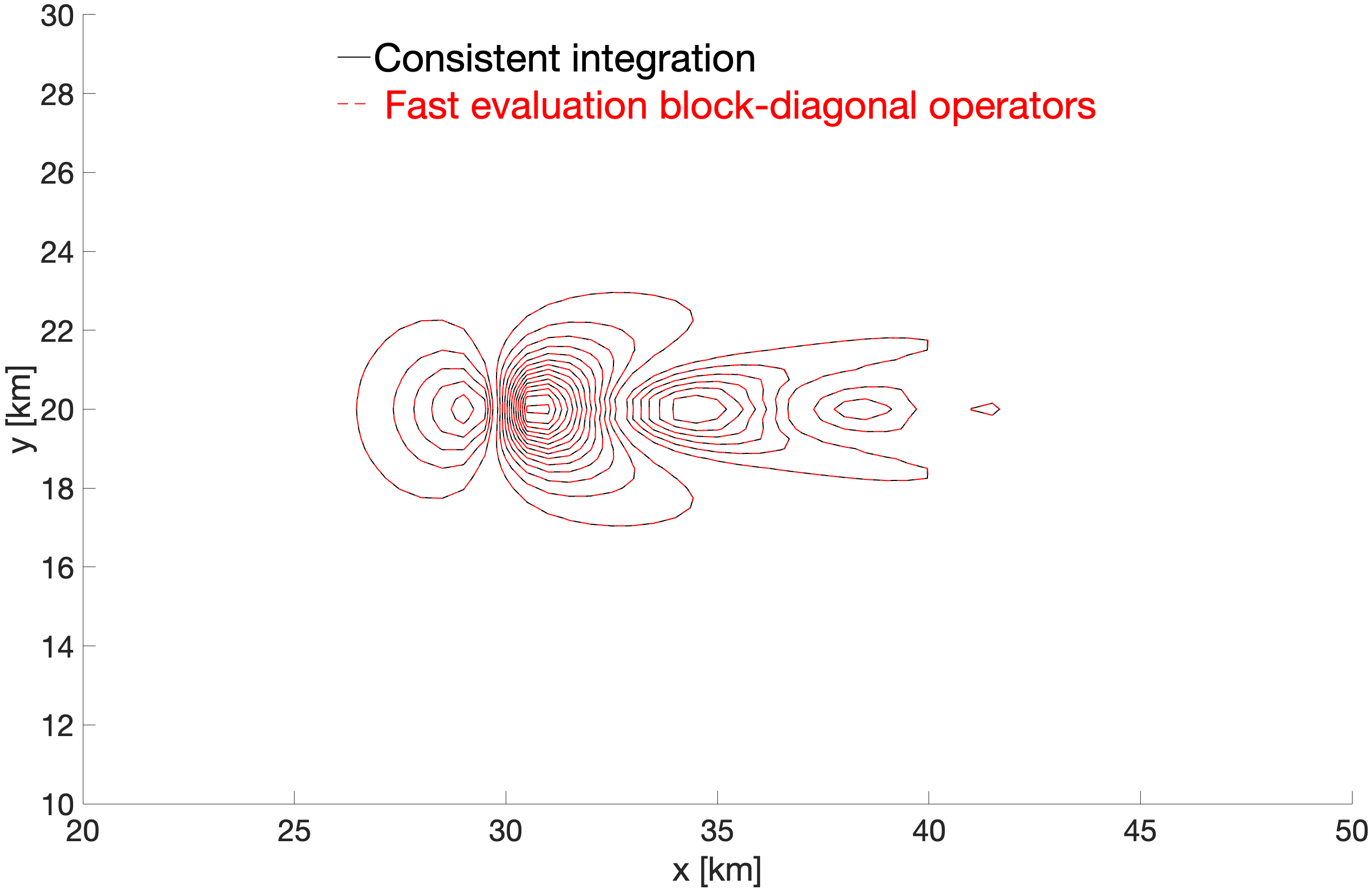}
    \end{subfigure}
    \caption{Flow over 3D hill test case, computed vertical velocity at $T_{f} = \SI{1}{\hour}$. Left: $x-z$ slice at $y = \SI{20}{\kilo\meter}$. Right: $x-y$ slice at $z = \SI{20}{\kilo\meter}$. Continuous black lines show the results with the consistent integration, whereas dashed red lines report the results with the fast evaluation of block-diagonal operators.}
    \label{fig:accuracy_block_diagonal}
\end{figure}

The strategy presented in \cite{kronbichler:2019} is already available in the \texttt{deal.II} library, but it was not considered in our previous work in order to avoid aliasing errors, as mentioned above. In this work we show the impact of this strategy on the numerical results and on the parallel performance. First, we analyze the impact of the use of the fast evaluation of block-diagonal operators in terms of accuracy. We consider a uniform mesh composed of $N_{el} = 60 \times 40 \times 16 = 38400$ elements, i.e. a resolution of \SI{250}{\meter} with a final time $T_{f} = \SI{1}{\hour}$. An excellent agreement is established between the two numerical strategies to solve the linear systems (Figure \ref{fig:accuracy_block_diagonal}). A computational time saving of around the $20\%$ is established using the fast evaluation of the block diagonal operators. Hence, this first test suggests that the fast evaluation of block-diagonal operators provides similar results in terms of accuracy in spite of the aliasing error. We aim to further investigate this matter in future work.

\begin{figure}[h!]
    \centering
    \includegraphics[width = 0.7\textwidth]{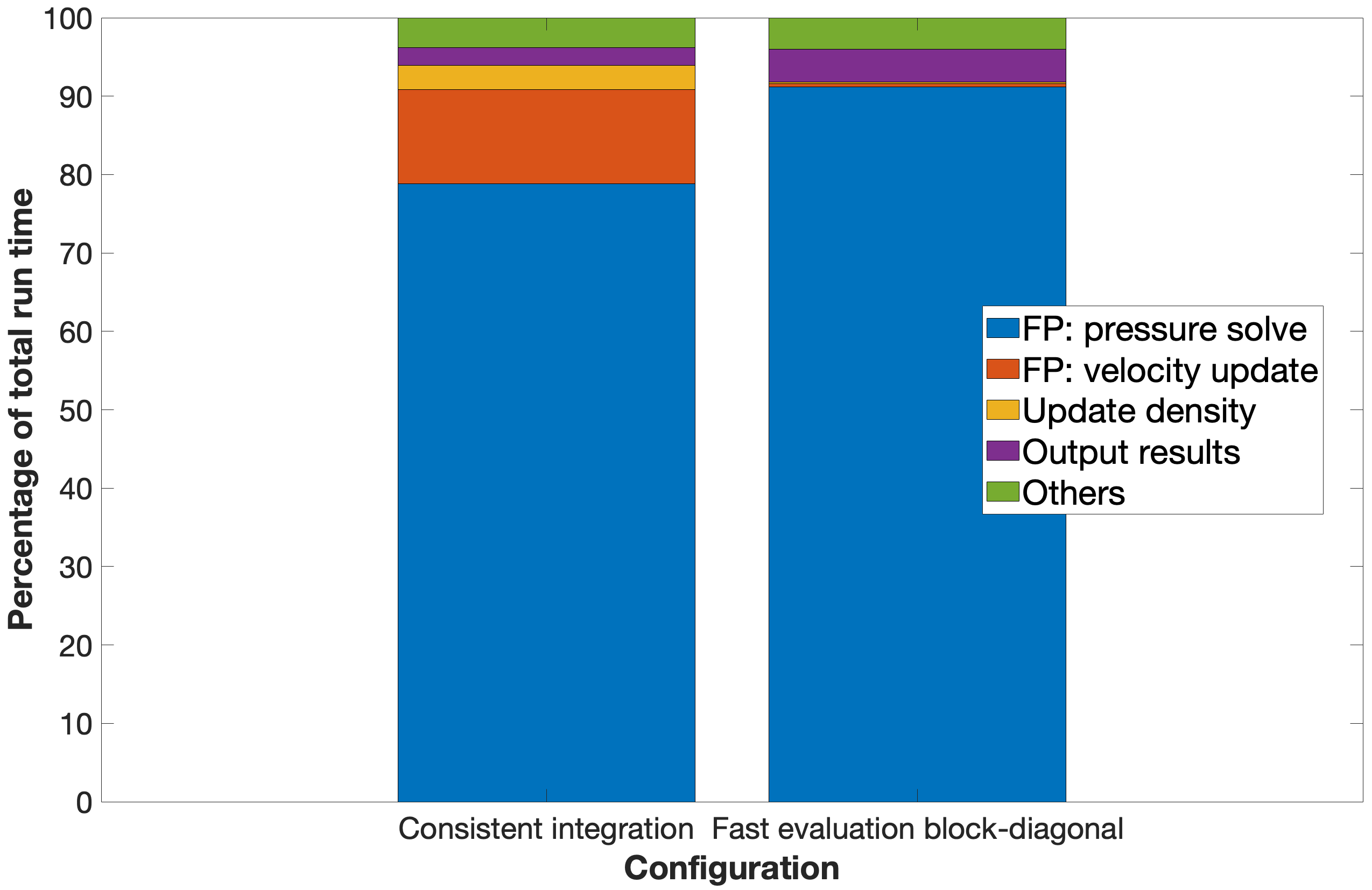}
    \caption{Distribution of the computational time spent in various blocks of the algorithm for the solution of problem \eqref{eq:stage_euler}, flow over a 3D hill test case with 2048 cores. Left: consistent numerical integration. Right: fast evaluation of block-diagonal operators.}
    \label{fig:profiling_block_diagonal}
\end{figure}

\begin{figure}[h!]
    \centering
    \begin{subfigure}{0.475\textwidth}
	\centering
        \includegraphics[width = 0.95\textwidth]{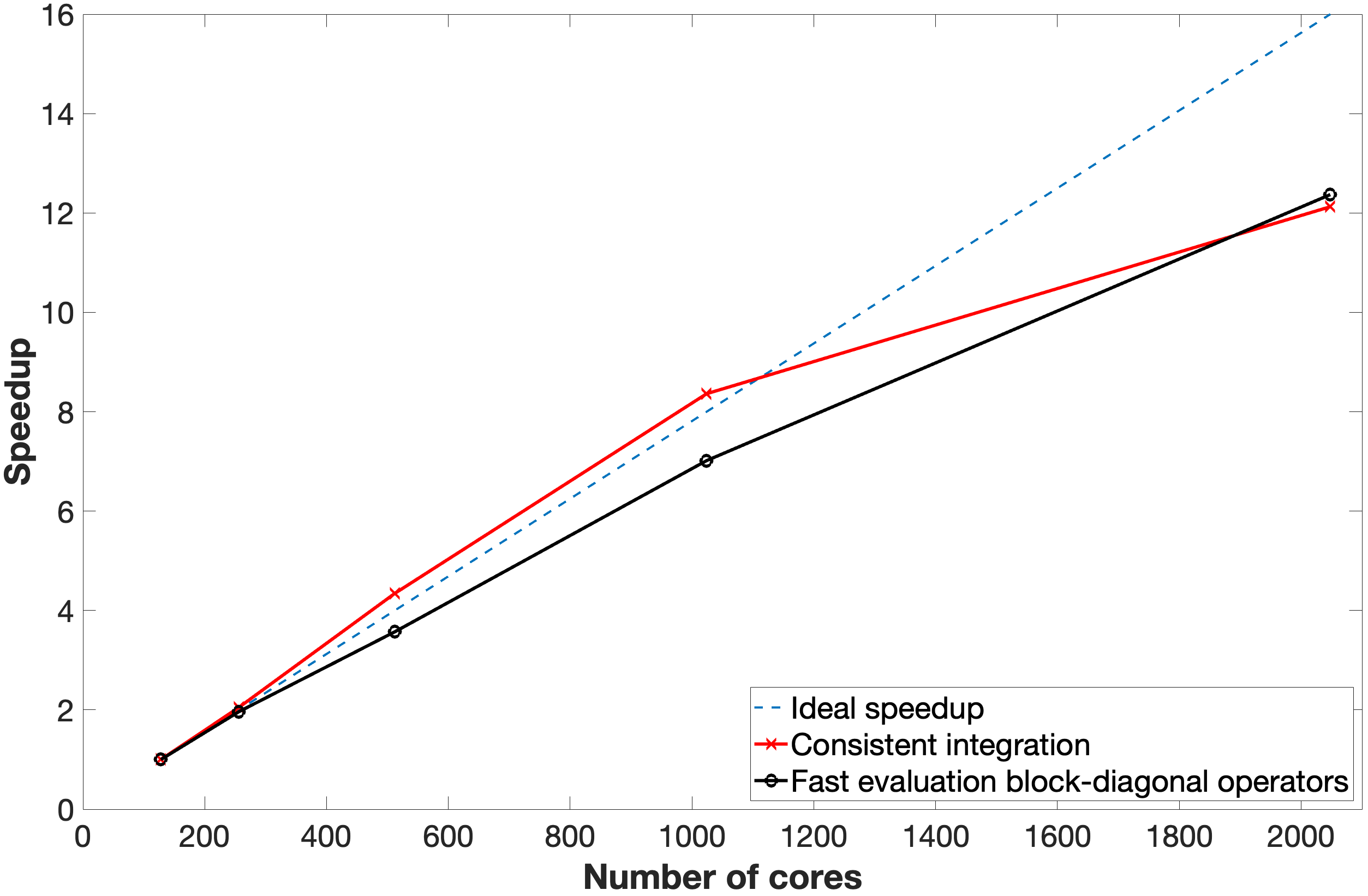}
    \end{subfigure}
    \begin{subfigure}{0.475\textwidth}
	\centering
        \includegraphics[width = 0.95\textwidth]{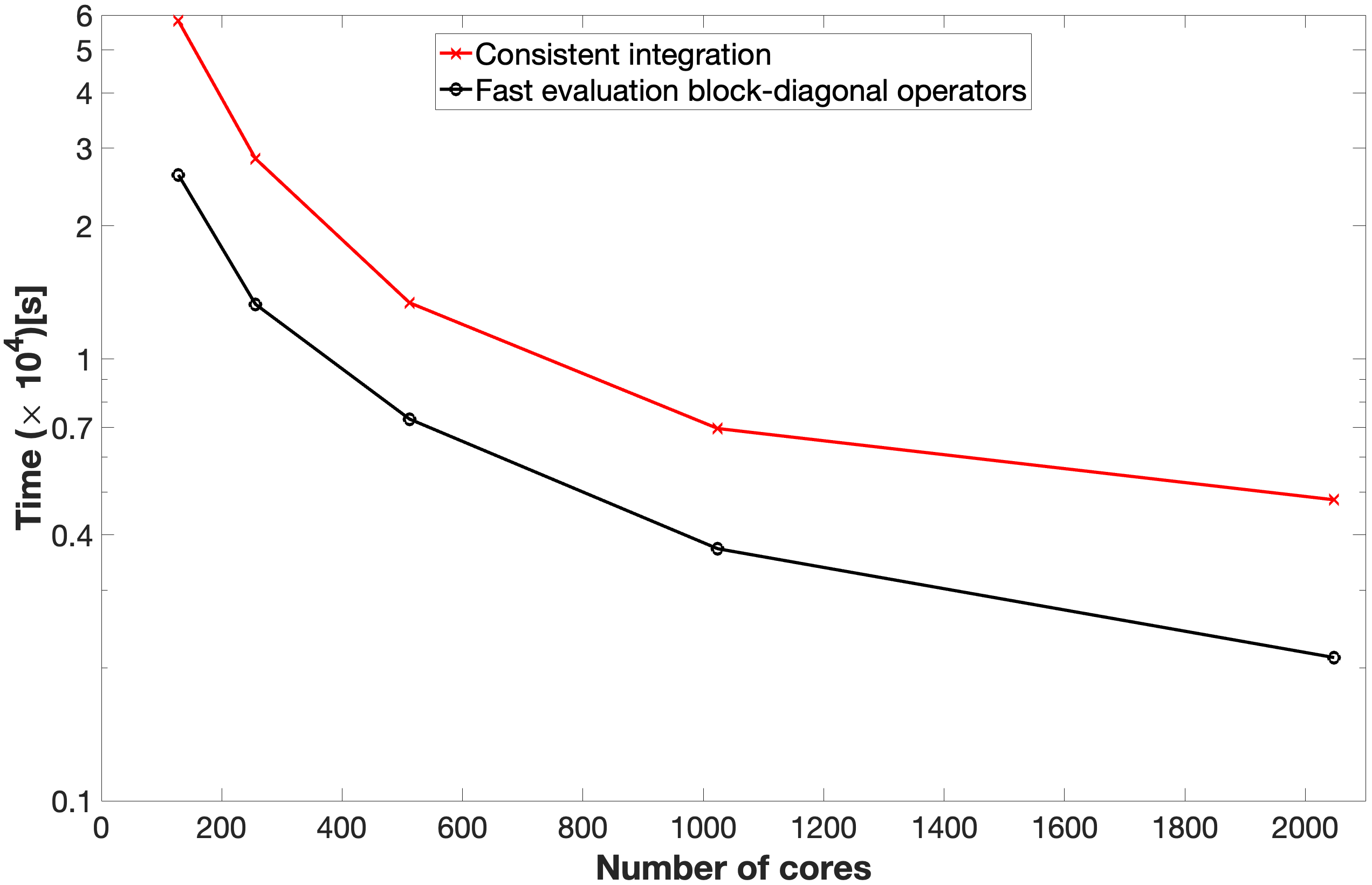}
    \end{subfigure}
    \begin{subfigure}{0.475\textwidth}
	\centering
        \includegraphics[width = 0.95\textwidth]{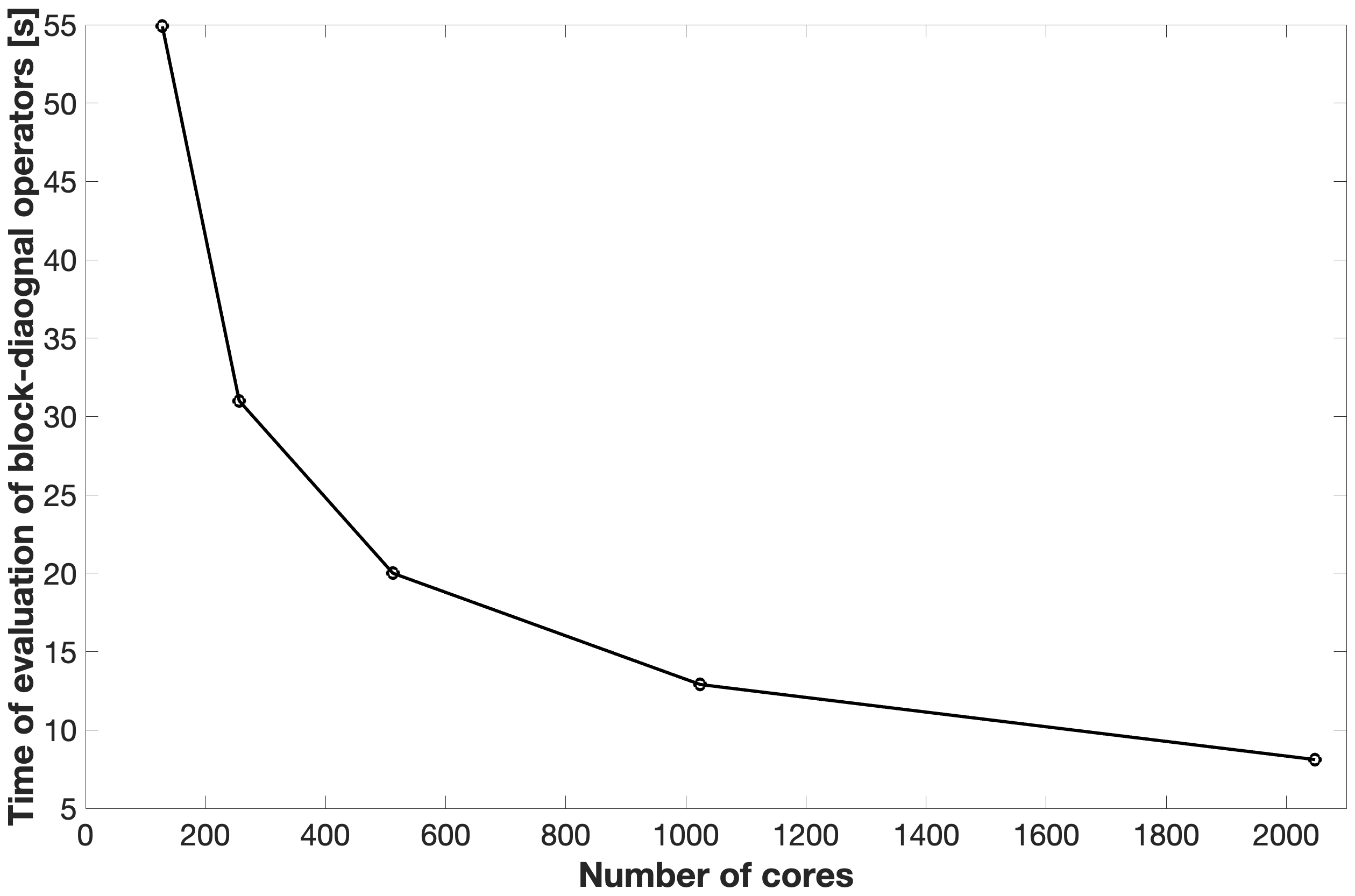}
    \end{subfigure}
    \begin{subfigure}{0.475\textwidth}
	\centering
        \includegraphics[width = 0.95\textwidth]{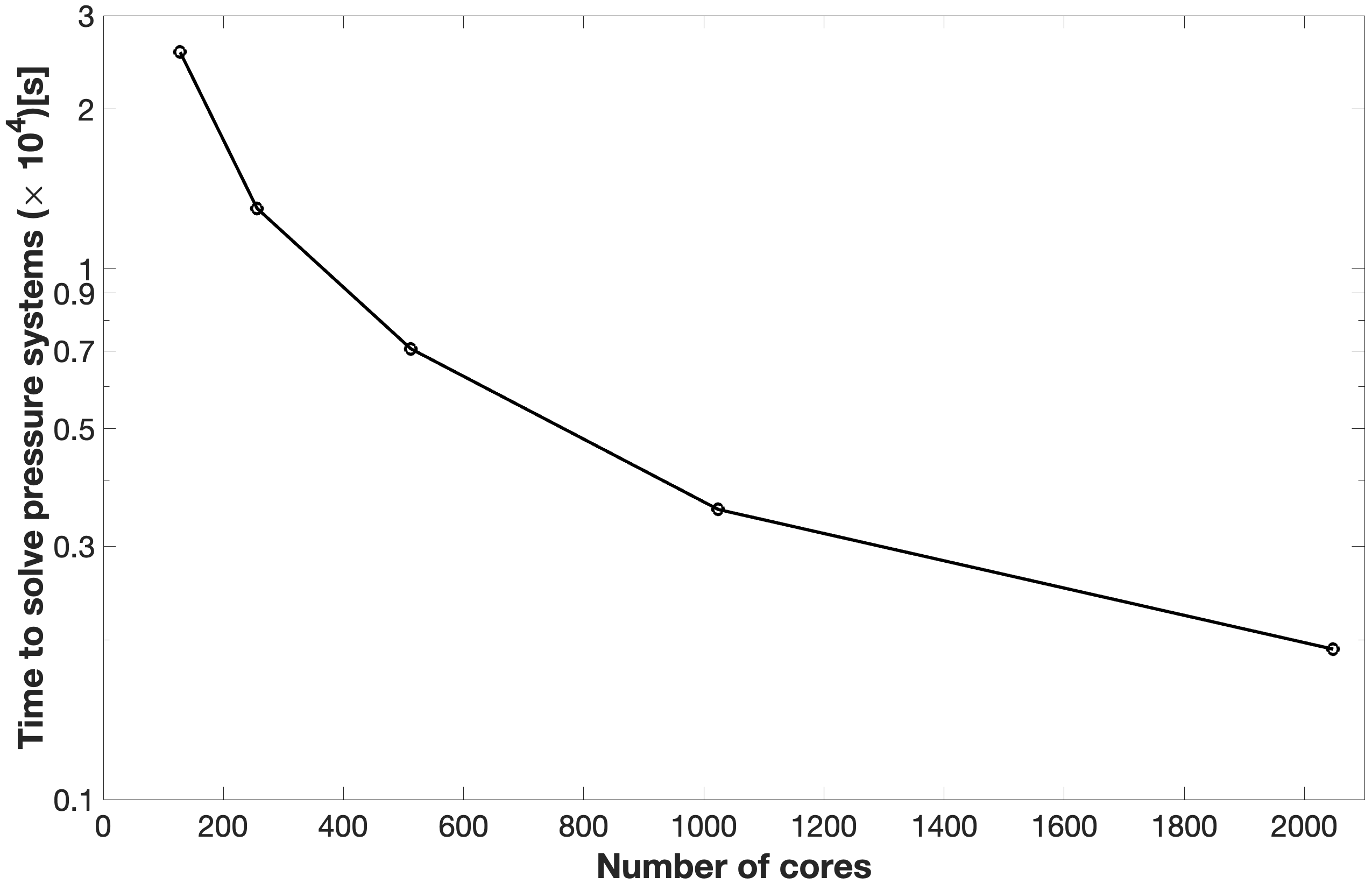}
    \end{subfigure}
    \caption{Parallel performance of the solver for problem \eqref{eq:stage_euler}, flow over a 3D hill test case with fast evaluation of block-diagonal operators. Top-left: strong scaling. Top-right: wall-clock times. Bottom-left: wall-clock times of evaluation of block-diagonal operators. Bottom-right: wall-clock times to solve linear systems for the pressure field \eqref{eq:fixed_point_discrete}. In the top panels, the black lines show the results with the fast evaluation of block-diagonal operators, while the red lines report the results obtained employing the consistent integration.}
    \label{fig:scaling_block_diagonal}
\end{figure}

Finally, we show the impact of the fast evaluation of the block-diagonal operators on the parallel performance. One can easily notice how the computational time spent in computing the velocity field is significantly reduced (orange bars in Figure \ref{fig:profiling_block_diagonal}). More specifically, the percentage with respect to the run time is reduced to around $0.52\%$ with 2048 cores. The strong scaling (Figure \ref{fig:scaling_block_diagonal}, top-left) when using the fast evaluation of the block-diagonal operators (red line) is apparently worse than when using consistent integration (black line) because the evaluation of the block-diagonal operators is so efficient that a small amount of time is required even for a relatively small number of cores (Figure \ref{fig:scaling_block_diagonal}, bottom-left). This consideration is further confirmed by the sizeable advantage in time-to-solution - $53\%$ with 128 cores, about $47\%$ with 2048 cores - brought by the fast evaluation of the block-diagonal operators (Figure \ref{fig:scaling_block_diagonal}, top-right). Moreover, one can easily notice that most of the computational time is required to solve \eqref{eq:fixed_point_discrete} to compute the the pressure field (Figure \ref{fig:scaling_block_diagonal}, bottom-right).

\section{Conclusions}
\label{sec:conclu}

This paper has reported the recent advances in the parallel performance of the IMEX-DG model for atmospheric dynamics simulations presented in \cite{orlando:2022, orlando:2023} and follows the analysis presented in \cite{orlando:2025b}. First, we have shown the impact of the use of a suitable preconditioner, then the impact in the use of suitable execution flags and latest version of \texttt{deal.II}, and, finally, the use of efficient matrix-free techniques for block-diagonal operators \cite{kronbichler:2019}.

In future work, we aim to further analyze the theoretical properties of the preconditioner, to further investigate the impact of the use of the direct inversion of the block-diagonal operator, and to exploit this technique to build a more suitable algebraic/geometric preconditioner since $(\mathbf{A}^{(l)})^{-1}$ can be readily evaluated. However, we point out that no theoretical result is available about the convergence properties of matrix-free multigrid preconditioners for non-symmetric linear systems arising from hyperbolic problems. Moreover, we plan to include more complex physical phenomena, in particular moist air, and to develop a three-dimensional dynamical core in spherical geometry including rotation so as to enable the testing on more realistic atmospheric flows. Finally, the data locality and high parallel efficiency of the DG method as well as the matrix-free approach which performs several computations on the fly makes the numerical method employed in this work particularly well-suited for a GPU implementation. Currently the \texttt{deal.II} library does not support discontinuous finite elements for GPU use. However, there is an ongoing effort of the \texttt{deal.II} community to provide a more complete and efficient GPU infrastructure including support for discontinuous finite elements. We aim to provide and test a GPU implementation once the library will be ready for this purpose.

\section*{Acknowledgements}

The simulations have been run thanks to the computational resources made available through the EuroHPC JU Benchmark And Development project EHPC-DEV-2024D10-054. We thank the Application Support EPICURE Team and in particular W.A. Mainassara for the support and help. This work has been partly supported by the ESCAPE-2 project, European Union’s Horizon 2020 Research and Innovation Programme (Grant Agreement No. 800897).

\printbibliography

@string{ QJRMS = "Quarterly Journal of the Royal Meteorological Society"}

@article{africa:2024,
         title     = {The deal.{II} library, {V}ersion 9.6},
         author    = {Africa, P.C. and Arndt, D. and Bangerth, W. and Blais, B. and Fehling, M. and Gassm{\"o}ller, R. and Heister, T. and Heltai, L. and Kinnewig, S. and Kronbichler, M. and Thiele, J.P and Turcksin, B. and Yushutin, V.},
         journal   = {Journal of Numerical Mathematics},
         volume    = {32},
         pages     = {369-380},
         year      = {2024},
         publisher = {De Gruyter}
}

@article{arndt:2023,
         title   = {The deal.{II} library, {V}ersion 9.5},
         author  = {Arndt, D. and Bangerth, W. and Bergbauer M. and Feder, M. and Fehling, M. and Heinz, J. and Heister, T. and Heltai, L. and Kronbichler, M. and Maier, M. and Munch, P. and Pelteret, J.-P. and Turcksin, B. and Wells, D. and Zampini, S.},
         pages   = {231-246},
         volume  = {31},
         journal = {Journal of Numerical Mathematics},
         year    = {2023}
}

@article{bangerth:2007,
         title   = {deal.{II}: a general-purpose object-oriented finite element library},
         author  = {Bangerth, W. and Hartmann, R. and Kanschat, G.},
         journal = {ACM Transactions on Mathematical Software },
         volume  = {33}, 
         pages   = {24-51},
         year    = {2007} 
}

@article{benacchio:2014,
         title     = {A blended soundproof-to-compressible numerical model for small-to mesoscale atmospheric dynamics},
         author    = {Benacchio, T. and O’Neill, W.P. and Klein, R.},
         journal   = {Monthly Weather Review},
         volume    = {142},
         pages     = {4416-4438},
         year      = {2014},
         publisher = {American Meteorological Society}
}

@book{giraldo:2020, 
      author    = {Giraldo, F.X.},  
      title     = {An {I}ntroduction to {E}lement-{B}ased {G}alerkin {M}ethods on {T}ensor-{P}roduct {B}ases.}, 
      chapter   = {}, 
      pages     = {}, 
      year      = {2020}, 
      publisher = {Springer {N}ature} 
}

@article{kennedy:2003,
         author  = {Kennedy, C.A. and Carpenter, M.H.},
         journal = {Applied Numerical Mathematics},
         pages   = {139-181},
         title   = {Additive {R}unge-{K}utta schemes for convection-diffusion-reaction equations},
         volume  = {44},
         year    = {2003}
}

@article{kronbichler:2016,
         title     = {Comparison of implicit and explicit hybridizable discontinuous {G}alerkin methods for the acoustic wave equation},
         author    = {Kronbichler, M. and Schoeder, S. and M{\"u}ller, C. and Wall, W.A.},
         journal   = {International Journal for Numerical Methods in Engineering},
         volume    = {106},
         pages     = {712-739},
         year      = {2016},
         publisher = {Wiley Online Library}
}

@article{kronbichler:2019,
         title     = {Fast matrix-free evaluation of discontinuous {G}alerkin finite element operators},
         author    = {Kronbichler, M. and Kormann, K.},
         journal   = {ACM Transactions on Mathematical Software (TOMS)},
         volume    = {45},
         pages     = {1-40},
         year      = {2019},
         publisher = {ACM New York, NY, USA}
}

@article{lanczos:1950,
         title   = {An iteration method for the solution of the eigenvalue problem of linear differential and integral operators},
         author  = {Lanczos, C.},
         journal = {Journal of research of the National Bureau of Standards},
         volume  = {45},
         pages   = {255-282},
         year    = {1950}
}

@article{lean:2024,
         title    = {The hectometric modelling challenge: {G}aps in the current state of the art and ways forward towards the implementation of 100-m scale weather and climate models},
         author   = {Lean, H.W. and Theeuwes, N.E. and Baldauf, M. and Barkmeijer, J. and Bessardon, G. and Blunn, L. and Bojarova, J. and Boutle, I.A. and Clark, P. A. and Demuzere, M. and others},
         journal = {Quarterly Journal of the Royal Meteorological Society},
         year    = {2024},
         volume  = {150},
         pages   = {4671-4708}
}

@article{melvin:2019,
         title   = {A mixed finite‐element, finite‐volume, semi‐implicit discretization for atmospheric dynamics: {C}artesian geometry},
         author  = {Melvin, T. and Benacchio, T. and Shipway, B. and Wood, N. and Thuburn, J. and Cotter, C.J.},
         journal = {Quarterly Journal of the Royal Meteorological Society},
         volume  = {145},
         pages   = {2835-2853},
         year    = {2019}
}

@article{orlando:2022,
         author  = {Orlando, G. and Barbante, P.F. and Bonaventura, L.},
         year    = {2022},
         title   = {An efficient {IMEX-DG} solver for the compressible {N}avier-{S}tokes equations for non-ideal gases},
         volume  = {471},
         pages   = {111653},
         journal = {Journal of Computational Physics}
}

@article{orlando:2023,
         title   = {An {IMEX-DG} solver for atmospheric dynamics simulations with adaptive mesh refinement},
         journal = {Journal of Computational and Applied Mathematics},
         pages   = {115124},
         year    = {2023},
         volume  = {427},
         author  = {Orlando, G. and Benacchio, T. and Bonaventura, L.}	 
}

@article{orlando:2024a,
         title   = {Robust and accurate simulations of flows over orography using non-conforming meshes}, 
         author  = {Orlando, G. and Benacchio, T. and Bonaventura, L.},
         year    = {2024},
         journal = QJRMS,
         volume  = {150},
         pages   = {4750-4770}
}

@article{orlando:2024b,
         title   = {Impact of curved elements for flows over orography with a {D}iscontinuous {G}alerkin scheme},
         journal = {Journal of Computational Physics},
         volume  = {519},
         pages   = {113445},
         year    = {2024},
         author  = {Orlando, G. and Benacchio, T. and Bonaventura, L.}
}

@article{orlando:2025a,
         title   = {Asymptotic-preserving {IMEX} schemes for the {E}uler equations of non-ideal gases},
         author  = {Orlando, G. and Bonaventura, L.},
         journal = {Journal of Computational Physics},
         volume  = {529},
         pages   = {113889},
         year    = {2025}
}

@article{orlando:2025b,
         title     = {Efficient and scalable atmospheric dynamics simulations using non-conforming meshes},
         author    = {Orlando, G. and Benacchio, T. and Bonaventura, L.},
         journal   = {Procedia Computer Science},
         volume    = {255},
         pages     = {33-42},
         year      = {2025},
         publisher = {Elsevier},
         note      = {Proceedings of the Second EuroHPC user day}
}

@misc{orlando:2025c,
      title         = {A quantitative comparison of high-order asymptotic-preserving and asymptotically-accurate {IMEX} methods for the {E}uler equations with non-ideal gases}, 
      author        = {Orlando, G. and Boscarino, S. and Russo, G.},
      year          = {2025},
      eprint        = {2501.12733},
      archivePrefix = {arXiv},
}

@article{pareschi:2005,
         title     = {Implicit-explicit {R}unge-{K}utta schemes and applications to hyperbolic systems with relaxation},
         author    = {Pareschi, L. and Russo, G.},
         journal   = {Journal of Scientific computing},
         volume    = {25},
         pages     = {129-155},
         year      = {2005},
         publisher = {Springer}
}

@article{steppeler:2003,
         author  = {Steppeler, J. and Hess, R. and Doms, G. and Sch{\"a}ttler, U. and Bonaventura, L.},
         title   = {Review of numerical methods for nonhydrostatic weather prediction models.},
         journal = {Meteorology and Atmospheric Physics},
         pages   = {287-301},
         year    = {2003},
         volume  = {82}
}

@book{varga:2009,
      author    = {Varga, R.S.},
      title     = {Matrix Iterative Analysis},
      year      = {2009},
      publisher = {Springer Science \& Business Media}
}

\end{document}